\documentclass[12pt]{amsart}
\usepackage{amsmath}
\usepackage{amsthm}
\usepackage{amssymb}
\usepackage{graphicx}
\usepackage{multirow}
\usepackage{lscape}

\newtheorem{rem}{Bemerkung}

\newcommand{\Q}{{\mathbb Q}}

\newcommand{\Z}{{\mathbb Z}}
\newcommand{\F}{{\mathbb F}}
\newcommand{\N}{{\mathbb N}}
\newcommand{\C}{{\mathbb C}}
\newcommand{\R}{{\mathbb R}}
\newcommand{\BP}{{\mathbb P}}

\newcommand{\La}{{\Lambda}}
\newcommand{\De}{{\Delta}}

\newcommand{\OM}{{\Omega}}
\newcommand{\Ga}{{\Gamma}}
\newcommand{\ga}{{\gamma}}
\newcommand{\bo}{{\boldsymbol{\omega}}}
\newcommand{\bzeta}{{\boldsymbol{\zeta}}}

\newcommand{\bF}{{\boldsymbol{F}}}

\newcommand{\ba}{{\boldsymbol{a}}}
\newcommand{\om}{{\omega}}
\newcommand{\bep}{{\boldsymbol{\epsilon}}}
\newcommand{\bof}{{\boldsymbol{f}}}

 \newcommand{\bx}{{\boldsymbol{x}}}

\newcommand{\bz}{{\boldsymbol{z}}}

\newcommand{\lra}{{\longrightarrow}}

\makeindex
\numberwithin{equation}{section}

\begin{document}
\setcounter{section}{-1}
\parindent=0pt

\title[On Drinfeld modular forms of higher rank III]
{On Drinfeld modular forms of higher rank III: The analogue of  the $k/12$-formula}
 \author{Ernst-Ulrich Gekeler}
 \maketitle
 
 \begin{abstract}
Continuing the work of \cite{7} and \cite{8}, we derive an analogue of the classical "$k/12$-formula" for Drinfeld modular forms of 
rank $r \geq 2$. Here the vanishing order $\nu_{\omega}(f)$ of one modular form at some point $\omega$ of the complex upper
half-plane is replaced by the intersection multiplicity $\nu_{\bo}(f_1,\ldots,f_{r-1})$ of $r-1$ independent Drinfeld modular forms
at some point $\bo$ of the Drinfeld symmetric space $\OM^r$. We apply the formula to determine the common zeroes of
$r-1$ consecutive Eisenstein series $E_{q^{i}-1}$, where $n-r<i<n$ for some $n \geq r$.
 \vspace{0.3cm}\end{abstract}

 \section{Introduction}
 Let $0 \not=f$ be an elliptic modular form of weight $k$ for $\Gamma = {\rm SL}(2,\Z)$. For $\om \in H$, the complex upper
 half-plane, $\nu_{\om}(f)$ denotes its vanishing order at $\om$. Then (see e.g. \cite{10}):
  $$\underset{\om\in \Ga\setminus H}{\sum{}^{\ast}} \nu_{\om}(f) +\frac{1}{2}\nu_i(f) + \frac{1}{3} \nu_{\rho}(f)+
   \nu_{\infty}(f) = k/12,\leqno{(0.1)}$$
where $\sum^{\ast}$ denotes the sum over a system of representatives of the non-elliptic classes of $H$ modulo the action of
$\Ga$, $i$ and $\rho$ are the primitive 4-th resp. 3-rd roots of unity in $H$, and $\nu_{\infty}(f)$ is the vanishing order at
infinity. Note that the denominator 2 resp. 3 is the excess $\# (\Ga_{\om})/2$ of the size of the stabilizer group  $\Ga_{\om}$
($\om \in \{i,\rho\}$), and $12 = 4\times 6/{\rm gcd}(4,6)$, where 4 and 6 are the weights of the canonical generators (usually
labelled $g_2,g_3$) of the $\C$-algebra of modular forms for $\Ga$.
 \vspace{0.3cm}
 
Now replace the data $\Z,\Q,\R,\C,H,\Ga$ with $A = \F_q[T]$, $K = {\rm Quot}(A) = \F_q(T)$, $K_{\infty} = \F_q((1/T))$, 
$C_{\infty} =$ completed algebraic closure of $K_{\infty}$, $\OM^2 = C_{\infty} \setminus K_{\infty}$, ${\rm GL}(2,A)$ as e.g.
in \cite{4}. Then a formula similar to (0.1) holds for Drinfeld modular forms $f$ for ${\rm GL}(2,A)$ on $\OM^2$, viz \cite{4}
5.14:
 $$\underset{\om \in {\rm GL}(2,A)\setminus \OM^2}{\sum{}^*} \nu_{\om}(f) + \frac{1}{q+1} \nu_{\epsilon}(f) + \frac{1}{q-1} 
  \nu_{\infty}(f) = k/(q^2-1), \leqno{(0.2)}$$
  where $\sum{}^*$ is the sum over the non-elliptic classes in $\OM^2$, $\epsilon$ is one fixed element of $\F_{q^2} \setminus \F_q$, and
  $\nu_{\infty}$ is the (appropriately defined) vanishing order of $f$ at $\infty$. The co-factor $1/(q-1)$ of $\nu_{\infty}$ comes from
  the facts that the Drinfeld discriminant $\Delta$ has a $(q-1)$-th root $h$ in the algebra of modular forms, and that the modular
  group is ${\rm GL}(2,A)$ instead of ${\rm SL}(2,A)$, see Remark 2.20 for more details. 
   \vspace{0.3cm}
   
The aim of the present work is to generalize (0.2) to Drinfeld modular forms of rank $r \geq 2$ as in \cite{7}  and \cite{8}.
Here the vanishing order $\nu_{\om}(f)$ of one modular form $f$ is replaced with the intersection number $\nu_{\bo}(\bof)$
of a system $\bof = \{f_1,\ldots,f_{r-1}\}$ of $r-1$ independent modular forms at $\bo \in \OM^r$, the $r-1$-dimensional Drinfeld
symmetric space. The result is the formula in Theorem 2.19. We give a handful of examples to show its usefulness. Among others, we
determine the common zeroes of the $r-1$ Eisenstein series $E_{q-1}$, where $n-r<i<n$ for some $n \geq r$.
 \vspace{0.5cm}
 
  {\bf Notation.} The notation largely agrees with that of \cite{7} and \cite{8}:
   \vspace{0.3cm}
   
$\F = \F_q =$ finite field with $q$ elements, of characteristic $p$;\\
$\overline{\F}  =$ algebraic closure of $\F$, $\F^{(n)} = \{x \in \overline{\F} ~|~x^{q^n} = x\}$;\\
$A=\F[T]$, $K = {\rm Quot}(A) = \F(T)$;\\
$C_{\infty}$ the completed algebraic closure of $K_{\infty} = \F((1/T))$;\\
$r \geq 2$ a natural number;\\
$\OM^r = \{\bo = (\om_1:\cdots:\om_r) \in \mathbb P^{r-1}(C_{\infty})~|~\mbox{the $\om_i$ are $K_{\infty}$-linearly indenpendent}\}$;\\
$\OM^r(R) = \{\bo \in \OM^r~|~\mbox{the $\om_i$ lie in the subring $R$ of $C_{\infty}$}\}$;\\
$\Ga = {\rm GL}(r,A)$, which acts on $\OM^r$ through fractional linear transformations, with center $Z \cong\F{}^*$;\\
${\rm Mod}(\Ga) = C_{\infty}[g_1,\ldots,g_r]$ (resp. ${\rm Mod}'(\Ga) = C_{\infty}[g_1,\ldots,g_{r-1},h]$) the ring of modular
forms of type zero (resp. of modular forms of arbitrary type).
 \vspace{0.3cm}
 
For any undefined notation and further explanation, we refer to \cite{7} and \cite{8}.
 
  \section{The moduli scheme $\overline{M}^r$}
  A Drinfeld $A$-module $\phi$ of rank $r \geq 2$ over $C_{\infty}$ is given by an operator polynomial
   $$\phi_T(X) = TX +g_1X^q+\cdots + g_rX^{q^r}\leqno{(1.1)}$$
with $g_i \in C_{\infty}$, $\Delta:= g_r \not=0$. Regarding the $g_i$ as indeterminates of weights $w_i:=q^{i}-1$, the moduli scheme of
such Drinfeld modules is
  $$M^r = {\rm Proj}(C_{\infty}[g_1,\ldots,g_r])_{g_r \not= 0},\leqno{(1.2)}$$
the complement of the zero locus $\partial M^r$ of $g_r$ in the weighted projective space  
 $$\overline{M}^r = {\rm Proj}({\rm Mod}(\Ga)),$$
where ${\rm Mod}(\Ga) = C_{\infty}[g_1,\ldots,g_r]$ is the graded ring of modular forms of type zero for $\Ga$. We also put
 $$M^1 = \overline{M}^1 = {\rm Proj}(C_{\infty}[g_1]) = \{{\rm point}\}.$$
 Then $\partial M^r$ is canonically isomorphic with $\overline{M}^{r-1}$, which allows induction procedures.
  \vspace{0.3cm}
  
(1.3) The Drinfeld symmetric space $\OM^r$ comes with a natural structure of rigid analytic space. With each $\bo \in \OM^r$ we can 
associate a homothety class of $A$-lattices $\La_{\bo} = A\om_1 + \cdots +A\om_r$ in $C_{\infty}$ and its attached class of Drinfeld
modules $\phi^{\bo}$ with
 $$  \phi_T^{\bo} = TX+g_1(\bo)X^q+\cdots+g_r(\bo)X^{q^r}.$$
 This way, the $g_i$ become modular forms of weights $w_i = q^{i}-1$ and type 0, and the map
  $$\begin{array}{llll}
  \beta: & \OM^r&\lra &M^r(C_{\infty})\\
   & \bo &\longmapsto& (g_1(\bo): \cdots:g_r(\bo))
  \end{array}$$
defines an isomorphism of the quotient analytic space $\Ga\setminus \OM^r$ with (the analytic space associated with)
$M^r$. Here and in the sequel, we make no difference in notation between algebraic data (varieties $V$ defined over
$C_{\infty}$) and analytic data (the corresponding analytic spaces $V^{an}$), or the underlying point sets.
 \vspace{0.3cm}  
 
(1.4) Let $\mathbb P$ be the projective space $\BP^{r-1}/C_{\infty} = {\rm Proj}(C_{\infty}[X_1,\ldots,X_r]) $ with coordinate functions
$X_1,\ldots,X_r$, and $\pi$ the map
 $$\begin{array}{rll}
\BP &\lra&\overline{M}^r\\
(x_1:\cdots:x_r) & \longmapsto& (x_1^{w_1}: \cdots :x_r^{w_r}).
 \end{array}$$
For a natural number $n$ coprime with the characteristic $p$, let $\mu_n$ be the group of $n$-th roots of unity in $C_{\infty}$.
Put
 $$H:= (\mu_{w_1} \times \cdots \times \mu_{w_r})/\mu_{q-1},$$
 where the group $\mu_{q-1}$ is diagonally embedded into $\prod \mu_{w_i}$. Then $H$ acts faithfully on $\BP$ by
 $(\ldots \epsilon_i \ldots)(\ldots x_i \ldots) = (\ldots \epsilon_i x_i \ldots)$, where $\epsilon_i \in \mu_{w_i}$, and $\pi$ is in fact
 the quotient map onto $H \setminus\BP = \overline{M}^r$.
  \vspace{0.3cm}
  
(1.5) Recall the elementary fact that for natural numbers $i$ and their associated $w_i = q^{i}-1$:
 $$i|j \Leftrightarrow w_i|w_j,$$
and  therefore for any non-empty $S \subset \N$:
 $$\gcd(w_i~|~i\in S) = w_s, \quad \mbox{where } s=\gcd(S).$$
 (1.6) For $\bx = (x_1:\cdots:x_r)\in \BP$ we let $S_{\bx} := \{i~|~1 \leq i \leq r,\, x_i \not= 0\}.$
 The stabilizer $H_{\bx}$ of $\bx$ in $H$ is
  $$H_{\bx} = \{\bep = (\epsilon_i) \in \prod_{1\leq i \leq r} \mu_{w_i}~|~i \in S_{\bx} \Rightarrow \epsilon_i=1\} \cdot 
  D_{\bx} /\mu_{q-1},$$
  where $D_{\bx}$ is the maximal constant diagonal subgroup of $\underset{i\in S_{\bx}}{\prod} \mu_{w_i}$. By (1.5), $D_{\bx}$ is the
  diagonally embedded group $\mu_{w_s} \hookrightarrow \underset{i \in S_{\bx}}{\prod} \mu_{w_i}$, $s = \gcd(S_{\bx})$, and so
   $$\#(H_{\bx}) = \underset{1\leq i \leq r \atop i \not\in S_{\bx}}{\prod} w_i \cdot w_s/w_1.$$
 In particular, $\bx$ is ramified under $\pi$ (i.e., $H_{\bx} \not= \{1\}$) if and only if at least one of $x_1,\ldots,x_r$ vanishes.
  \vspace{0.3cm}  
 
(1.7) A modular form of weight $k$ and type zero for $\Ga$ is an isobaric polynomial of weight $k$ in $g_1,\ldots,g_r$, see \cite{8} 1.8.
We regard it as a form of weight $k$ on $\overline{M}^r$, that is, a function $f$ on $C_{\infty}^r \setminus \{0\}$
that satisfies $f(c^{w_1}g_1,\ldots,c^{w_r}g_r) = c^k f(g_1,\ldots,g_r)$ for $c \in C_{\infty}^*$. For such an $f$, the pullback
$\pi^*f$ is a homogeneous polynomial of degree $k$ in $X_1,\ldots,X_r$, and the automorphy condition translates to
$(\pi^*f)\circ \gamma = \pi^*f$ for $\gamma \in \Ga$, which acts on $\BP$ through left matrix multiplication.
 \vspace{0.3cm}
 
 (1.8) We adapt these considerations to modular forms of non-zero type. There exists a function $h$ of $\OM^r$ such that
 $h^{q-1} = (-1)^r/T \cdot \De$, see (\cite{7}, \cite{1}), and $h$ is a modular form of weight $w'_r:= (q^r-1)/(q-1)$ and type 1.
 For simplicity, we scale $h$ such that $h^{q-1} = \De = g_r$ (i.e., we replace $h$ with $((-1)^rT)^{1/(q-1)}h$, and abuse the
 name $h$ for the resulting form). Now the algebra ${\rm Mod}'(\Ga)$ of all modular forms (regardless of type) for $\Ga$ is the graded
 algebra $C_{\infty}[g_1,\ldots,g_{r-1},h]$ (\cite{8} 1.8), which suggests to also consider the weighted projective space
  $$\overline{M}' = \overline{M}'{}^r := {\rm Proj}({\rm Mod}'(\Ga)),$$
where the weight of the indeterminate $h$ is $w'_r = (q^r-1)/(q-1)$. Henceforth, we suppress the superscript $r$ if no
ambiguity arises. Let $M' = {M'}^r$ be the complement in $\overline{M}'$ of the divisor ($h=0$). Similar to (1.3), the map
 $$\begin{array}{rll}
  \beta' :\: \OM = \OM^r &\lra & M'(C_{\infty})\\
   & \longmapsto & (g_1(\bo): \cdots : g_{r-1}(\bo):h(\bo))
 \end{array}$$
 identifies $\Ga'\setminus \OM$ with $M'(C_{\infty})$, where $\Ga' := {\rm SL}(r,A) \subset \Ga$. Apparently, $\pi':\: \BP \lra \overline{M}'$
 defined by $(x_1:\cdots:x_r) \longmapsto (x_1^{w_1}:\cdots: x_r^{w'_r})$ is the Galois subcover of $\pi:\: \BP \lra \overline{M}$ that
 corresponds to the image $H'$ of $\underset{1\leq i \leq r-1}{\prod} \mu_{w_1} \times \mu_{w'_r}$ in $H =
 {\rm Gal}(\BP|\overline{M})$. We have 
  $$\deg (\overline{M}'|\overline{M}) = [\Ga:\Ga' \cdot Z] = [H:H'] = \gcd(q-1,r). \leqno{(1.9)}$$
  Similar to (1.7), a modular form of weight $k$ and type $m \in \Z/(q-1)$ is a form of weight $k$ on $\overline{M}'$. It may be written as
  an isobaric polynomial of weight $k$ in $g_1,\ldots,g_{r-1},h$ of shape $h^m$ times a modular form of weight $k-m \circ w'_r$ and type 0,
  where $0 \leq m < q-1$ is the standard representative. The pullback ${\pi'}^*f$ is a homogeneous polynomial in 
  $X_1,\ldots,X_r$ and satisfies 
   $$({\pi'}^*f)\circ \ga = (\det \ga)^m({\pi'}^* f)$$
   for $\ga \in \Ga$.
   
    \section{The "$k/12$-formula".}
    The number $r \geq 2$ is fixed throughout, and $M = M^r$, $\overline{M} = \overline{M}^r$, $\OM = \OM^r$, etc.
 \vspace{0.3cm}    
 
 (2.1) We call a family $\bof = \{f_1,\ldots,f_{r-1}\}$ of $r-1$ modular forms for $\Ga$ {\em independent} if their common vanishing
 locus $\underset{1\leq i \leq r-1}{\bigcap} V(f_i) \subset \overline{M}$ has dimension zero. A common zero of 
 $f_1\ldots,f_{r-1}$ is briefly called a {\em zero of } $\bof$.
  \vspace{0.3cm}
  
In what follows, let $R$ be a noetherian local $C_{\infty}$-algebra with maximal ideal $\mathfrak m$ and residue class field
$C_{\infty}$, and regular of dimension $n \geq 1$. In practice, $R$ will be the local ring $\mathcal O_{\bo,\OM}$ for some
$\bo \in \OM$, or the (algebraic or analytic) local ring $\mathcal O_{\bx,\BP}$ of $\BP = \BP^{r-1}$ at $\bx$. In these cases,
$n = r-1$. 
For a system $\bof = \{f_1,\ldots,f_n\}$ of parameters of $R$ (i.e., for large $k$, $\mathfrak m^k$ is contained in the ideal
$\langle \bof\rangle_R$ generated by $\bof$) we define the {\em intersection number of $\bof$ at $R$} as the length of the
Artin ring $R/\langle \bof \rangle_R$, which in the given case is the dimension 
 $$\nu(R,\bof) := \dim_{C_{\infty}} R/\langle \bof \rangle_R. \leqno{(2.2)}$$
 Some of the relevant properties are:
  \vspace{0.3cm}
  
(2.3) The completion $\hat{R}$ of $R$ shares the requirements on $R$ (noetherian, local, regular of dimension $n$) and 
   $$\nu(\hat{R},\bof) = \nu(R,\bof);$$
(2.4) If $R$ is the local ring of a $C_{\infty}$-surface (i.e., $n=2$) at a smooth point, $\nu(R,\bof)$ agrees with the naive intersection
number as defined in \cite{9} V Sect.\,1;
    \vspace{0.3cm}
    
(2.5) If $R'$ is a subring subject to the same requirements as $R$ and such that $R|R'$ is finite and flat (hence free of some rank
$e \in \N$), and $\bof$ is a system of parameters for $R'$, then it remains so for $R$, and the formula
 $$\nu(R,\bof) = e \cdot \nu(R',\bof)$$
 holds;    
 \vspace{0.3cm}
 
(2.6) If $\bof = \{f_1,\ldots,f_n\}$ and $\bof' = \{f'_1,f_2,\ldots,f_n\}$ are such that $\nu(R,\bof)$ and $\nu(R,\bof')$ are defined,
so is $\nu(R,\bof'')$ with $\bof'' = \{f_1f'_1,f_2,\ldots,f_n\}$, and $\nu(R,\bof'') = \nu(R,\bof)+\nu(R,\bof').$ 
 \vspace{0.3cm}
 
Here (2.3) is the synopsis of several well-known facts that can be found e.g. in Bourbaki, Alg\`{e}bre commutative \cite{2} Ch. VIII,
Sect. 5 or \cite{3} Sect. 7, 10, 19, (2.4) is \cite{9} I 5.4, and (2.5) is obvious, as
$R/\langle \bof\rangle_R \cong (R'/\langle \bof \rangle_{R'})^{e}$ as an $R'$-module. For (2.6), consider the canonical exact sequence
of finite-dimensional $C_{\infty}$-spaces 
 $$0 \lra \langle \bof \rangle/\langle \bof'' \rangle_R \lra R/\langle \bof'' \rangle_R \lra R/\langle\bof \rangle_R \lra 0$$
 and the map
  $$\varphi:\: R/\langle \bof' \rangle_R \lra \langle \bof \rangle _R/\langle \bof''\rangle_R$$
induced from $a \longmapsto af_1$ ($a\in R$). Certainly, $\varphi$ is well-defined and surjective, and its injectivity is a
consequence of unique factorization in the regular ring $R$ (\cite{3}, Theorem 19.19). Now the result follows from comparing
dimensions.
 \vspace{0.3cm}
 
{\bf 2.7 Remarks.} (i) As the algebraic local ring $\mathcal O_{\bx,V}$ of a $C_{\infty}$-variety $V$ and its ''analytification''
$\mathcal O_{\bf x,V^{an}}$ share their  completions, we need by virtue of (2.3) not distinguish between algebraic and analytic 
intersection numbers $\nu(.,\bof)$.\\
(ii) The conditions of (2.5) are fulfilled if $R'$ is the subring of invariants (required to be regular) of a finite group $H$ of
automorphisms of $R$. In this case, $e={\rm rank}_{R'}(R) = \#(H)$.
 \vspace{0.3cm}
 
 (2.8) Let now a family $\bof$ of $r-1$ independent modular forms $f_1,\ldots,f_{r-1}$ with weights $k_1,\ldots,k_{r-1}$ be given, regarded as
 forms  on $\overline{M}'$ (see (1.8)). We write $f_i^*$ for the pullback ${\pi'}^*(f_i)$ to $\BP$, a homogeneous element of weight
 $k_i$ of $C_{\infty}[X_1,\ldots,X_r]$, and $\bof^*$ for the family of the $f_i^*$. From the generalized B\'{e}zout theorem 
 (\cite{9} I, 7.7) we find
  $$\underset{\bx\in \BP(C_{\infty})}{\sum}\nu_{\bx}(\bof) = \underset{1 \leq i \leq r-1}{\prod}k_i,\leqno{(2.8.1)}$$
where $\bx$ runs through the points of $\BP = \BP^{r-1}$ and $\nu_{\bx}(\bof)$ is the intersection number
$\nu(\mathcal O_{\bx,\BP},\bof^*)$, which is positive if and only if the $f_i^*$ vanish simultaneously in $\bx$. (As usual,
a "form" $f$ is made a "function" $f^{(i)}$ around $\bx$ through dehomogenization w.r.t. to the 
$i$-th coordinate, where $x_i \not= 0$, and $\nu(\cdot,\bof^*) = \nu(\cdot,\bof^{*(i)})$ doesn't depend on the choice of 
$i$.)
 \vspace{0.3cm}
 
Recall that $\pi:\: \BP \lra \overline{M} = H \setminus \BP$ is the quotient map by the group $H$ of (1.4). Let $\eta$ be an element of
$H$. Then 
 $$\nu_{\eta\bx}(\bof) = \nu(\mathcal O_{\eta\bx},\bof^*) = \nu(\mathcal O_{\bx},\bof^* \circ \eta^{-1}) =
  \nu(\mathcal O_{\bx},\bof^*)= \nu_{\bx}(\bof),$$
as $\eta^{-1}$ transforms $f_i^*$ to a constant multiple of $f_i^*$. That is, $\nu_{\bx}(\bof)$ depends only on the
$H$-orbit of $\bx$, i.e., on the point $\bz := \pi(\bx) \in \overline{M}(C_{\infty})$. Putting
 $$\nu_{\bz}(\bof) := \# (H_{\bx})^{-1} \nu_{\bx}(\bof)\leqno{(2.8.2)}$$
 with the stabilizer $H_{\bx}$ of $\bx$, we find
  $$\underset{\bx' \in H\bx}{\sum} \nu_{\bx'} (\bof) = \# (H) \nu_{\bz}(\bof) \leqno{(2.8.3)}$$
for each fixed $\bx \in \BP(C_{\infty})$. Summing over the $H$-orbits in $\BP(C_{\infty})$ and dividing by $\#(H)$, we find as
a first approach to the wanted formula the following
 $$\underset{\bz \in \overline{M}(C_{\infty})}{\sum} \nu_{\bz}(\bof) = \#(H)^{-1} \underset{1 \leq i \leq r-1}{\prod} k_i\,.
  \leqno{(2.9)}$$
{\bf Remark.} Although the definition (2.8.2) of $\nu_{\bz}(\bof)$ is motivated from (2.5) and the resulting (2.8.3), it 
wouldn't make sense to define $\nu_{\bz}(\bof)$ as $\nu(\mathcal O_{\bz,\overline{M}},\bof)$, since, due to the presence 
of weights, the $f_i$ cannot be regarded as functions on $\overline{M}$ around $\bz$. In particular, $\nu_{\bz}(\bof)$ as defined
by (2.8.2) is in general not an integer, see Remark 2.20.
 \vspace{0.3cm}  
 
(2.11) From now on, we assume the homogeneous coordinates \mbox{$(\om_1: \cdots:\om_r)$} on $\OM$ normalized by $\om_r=1$. Then
the $f_i$ become functions on $\OM$. Suppose that $\bz := \pi(\bx)$ belongs to $M(C_{\infty})$, image  of
$\bo \in \OM$. What is the relationship between $\nu_{\bo}(\bof) := \nu(\mathcal O_{\bo,\OM},\bof)$ and 
$\nu_{\bx}(\bof)$ or $\nu_{\bz}(\bof)$?
 \vspace{0.3cm}
 
In order to answer the question, we construct a map $\alpha$ from a small neighborhood $U$ of $\bx^{(0)} \in \BP(C_{\infty})$ to
a neighborhood of $\bo^{(0)} \in \OM$ such that $\alpha(\bx^{(0)})=\bo^{(0)}$ and the diagram
 $$ \begin{array}{ccc}
  &U &\\
  \alpha\swarrow & & \searrow \pi\\
  \OM & \stackrel{\beta}{\lra}& M
  \end{array}\leqno{(2.11.1)}$$
commutes. The wanted relationship will come out by looking at the structure of the $\mathcal O_{\bo^{(0)},\OM}$-algebra 
 $\mathcal O_{\bx^{(0)},\BP}$. 
  \vspace{0.3cm}
  
As $h(\bo^{(0)}) \not=0$, there exists a small neighbourhood $V$ of $\bo^{(0)}$ and a holomorphic function $\tilde{h}$ on $V$ such
that $\tilde{h}^{w'_r} = h$. By \cite{8} Proposition 3.15 the $r-1$ functions $g_1,\ldots,g_{r-1}$ are local coordinates around
$\bo^{(0)}$, that is the $g_i-g_i(\bo^{(0)})$ with $1 \leq i \leq r-1$ form a regular sequence in the local ring
$\mathcal O_{\bo^{(0)},\OM}$. This remains true upon replacing the $g_i$ with $\tilde{g}_i := g_i\tilde{h}^{-w_i}$.
As $x_r^{(0)} \not=0$, we may consider a neighborhood $U$ of $\bx^{(0)}$ of shape
$U = \{\bx \in \BP(C_{\infty}) ~|~|\frac{x_i}{x_r}-\frac{x_i^{(0)}}{x_r^{(0)}}| \leq \epsilon,\: 1 \leq i \leq r-1\}$ for small
$\epsilon >0$. Choosing first $V$ and then $U$ sufficiently small, we may assume that 
 \begin{itemize}
  \item the $\tilde{g}_i$ are defined on $V$ and are coordinates on $V$;
  \item the map $\alpha:\: U \lra V$ defined below is in fact well-defined;
  \item for $\eta \in H$, $\eta(U) \cap U \not= \emptyset$ implies $\eta \in H_{\bx^{(0)}}$, and 
  \item $\eta(U) = U$ for $\eta \in H_{\bx^{(0)}}$. 
\end{itemize}
Here $\alpha:\: U \lra V$ maps $\bx = (x_1:\cdots:x_r) \in U$ to the point $\bo \in V$ with $\tilde{g}$-coordinates
$\tilde{g}_i(\bo)$, where
 $$\tilde{g}_i(\bo) = (\frac{x_i}{x_r})^{w_i} \quad (1 \leq i \leq r-1).\leqno{(2.11.2)}$$
By construction, for $\bo = \alpha(\bx)$
\begin{eqnarray*}
  \ga(\bo) &=& (g_1(\bo): \cdots:g_{r-1}(\bo):h^{q-1}(\bo))\\
   &=& (\tilde{g}_1(\bo):\cdots: \tilde{g}_{r-1}(\bo):1)\\
    &=& ((\frac{x_1}{x_r})^{w_1} : \cdots : (\frac{x_r}{x_r})^{w_r})\\
     &=& (x_1^{w_1}: \cdots: x_r^{w_r})\\
      &=& \pi(\bx),
 \end{eqnarray*}
 and (2.11.1) in fact commutes. Now $\pi|_U :\: U \lra \pi(U)$ and $\beta|_{\alpha(U)}:\:\alpha(U) \lra \pi(U)$ are faithfully flat as
 quotient morphisms, and so is $\alpha:\:U \lra \alpha(U)$. In particular, $R:= \mathcal O_{\bx^{(0)},\BP} = 
 \mathcal O_{\bx^{(0)},U}$ is finite and faithfully flat as an algebra over $R' := \mathcal O_{\bo^{(0)},\OM} = 
 \mathcal O_{\bo^{(0)},V}$, hence free of rank ${\rm rank}_{R'} (R)= \deg\,\alpha$.
  \vspace{0.3cm}
  
(2.12) We next determine the degrees of the maps $\alpha,\beta,\pi$. Let $\La_{\bo^{(0)}} = A\om_1^{(0)} + \cdots + 
A\om_r^{(0)}$ be the lattice determined by $\bo^{(0)}$, with associated Drinfeld module $\phi_T^{\bo^{(0)}}(X) = TX +
\underset{1 \leq i \leq r}{\sum} g_i(\bo^{(0)}) X^{q^{i}}$. Then
 \begin{eqnarray*}
 \Ga_{\bo^{(0)}} &\cong &\{c \in C_{\infty}^*~|~c\La_{\bo^{(0)}}=\La_{\bo^{(0)}}\} = {\rm Aut}(\phi^{\bo^{(0)}}) \\
   &= & \{c \in C_{\infty}^*~|~c^{w_i} g_i(\bo^{(0)}) = g_i(\bo^{(0)})\: \forall\: i\} \\
   &= &\mu_{w_s}.
   \end{eqnarray*}
Here $s = \gcd\{i~|~1 \leq i \leq r,\: g_i(\bo^{(0)}) \not=0\}$, where we have used (1.5). As $\Ga/Z$ acts faithfully on $\OM$, the
local ring $R'' := \mathcal O_{\bz^{(0)}),\overline{M}}$ at $\bz^{(0)} := \beta(\bo^{(0)})$ is the ring of invariants under
$\Ga_{\bo^{(0)}}/Z \cong \mu_{w_s}/\mu_{q-1}$ in $R'$, and ${\rm rank}_{R''}(R') = w_s/w_1 = (q^s-1)/(q-1)$.
 \vspace{0.3cm}
  
(2.13) We have determined $H_{\bx^{(0)}}$ in (1.6). This yields 
 $${\rm rank}_{R''}(R) = \#(H_{\bx^{(0)}}) = \underset{1\leq i \leq r \atop x_i^{(0)} \not= 0}{\prod} w_i \cdot w_s/w_1.$$  
(2.14) Let $\eta$ be an element of $H_{\bx^{(0)}}$, represented modulo $\mu_{q-1}$ through entries $\eta_i$, where
$\eta_i \in \mu_{w_i}$ if $x_i^{(0)} = 0$ and $\eta_i = c$ with some constant $w_s$-th root of unity $c$ if $x_i^{(0)} \not= 0$. Then
  \begin{eqnarray*}
 \alpha \circ \eta = \alpha &\Leftrightarrow &\forall\:\bx \in U,\: \forall \: i = 1,\ldots,r-1 (\frac{\eta_i}{\eta_r} \frac{x_i}{x_r})^{w_i}
 =(\frac{x_i}{x_r})^{w_i}\\
  & \Leftrightarrow & c \in \mu_{q-1}, \mbox{ as } \eta_r = c\\
   & \Leftrightarrow &  \eta \in (\underset{1\leq i \leq r \atop x_i^{(0)}=0}{\prod} \mu_{w_i} \times \mu_{q-1} )/\mu_{q-1} =: 
   H_{\bo^{(0)}}.
 \end{eqnarray*}
 Therefore, ${\rm rank}_{R'}(R) = \#(H_{\bo^{(0)}}) = \# (H_{\bx^{(0)}})/\#(\Ga_{\bo^{(0)}}/Z)$. 
  \vspace{0.3cm}
 
(2.15) By its very construction, the pullback $\alpha^*(\tilde{g}_i)$ of (the germ of) $\tilde{g}_i \in R' = \mathcal O_{\bo^{(0)},\OM}$
in $R = \mathcal O_{\bx,\BP}$ is $(X_i/X_r)^{w_i}$ ($1 \leq i \leq r-1$), see (2.11.3). If thus $f = F(g_1,\ldots,g_{r-1},h)$ is an 
arbitrary modular form with an isobaric $C_{\infty}$-polynomial  $F$, then $\alpha^*(f) = F\left( 
(\frac{X_1}{X_r})^{w_1},\ldots,(\frac{X_{r-1}}{X_r})^{w_r},1\right)$, and for an independent family $\bof = \{f_1,\ldots,f_{r-1}\}$ of modular
forms,
 $$\nu_{\bx^{(0)}}(\bof) = \#(H_{\bo^{(0)}}) \nu_{\bo^{(0)}}(\bof), \leqno{(2.16)}$$
 taking (2.5) and Remark 2.7\,(ii) into account.
  \vspace{0.3cm}
  
Having answered the question in (2.11), we may summarize what has been shown.
 \vspace{0.3cm}  
 
{\bf 2.17 Proposition.} {\it 
Let $\bo \in \OM$ and $\bx = (x_1:\cdots: x_r) \in \BP(C_{\infty})$ be such that $\pi(\bx) = \bz = \beta(\bo)$. Let further
$\bof = \{f_1,\ldots,f_{r-1}\}$ be an independent family of modular forms for $\Ga$. Then the intersection numbers
satisfy
 $$\frac{\nu_{\bo}(\bof)}{\#(\Ga_{\bo}/Z)} = \frac{\nu_{\bx}(\bof)}{\#(H_{\bx})} = \nu_{\bz}(\bof).$$
Here $\#(\Ga_{\bo}/Z) = w_s/w_1$ for some divisor $s$ of $r$ and $\#(H_{\bx})= (\underset{1\leq i \leq r \atop x_i=0}{\prod} w_i)
w_s/w_1$ with $w_i= q^{i}-1$.} 
 \begin{proof}
 The second equality is the definition of $\nu_{\bz}$, and the first one is (2.16) together with $H_{\bx}/H_{\bo} \cong \Ga_{\bo}/Z$.
 \end{proof}
 
{\bf 2.18 Remark.}  The map $\alpha$ used in (2.11)--(2.15) and therefore the $\mathcal O_{\bo,\OM}$-structure of 
$\mathcal O_{\bx,\BP}$ depends on the choice of the root $\tilde{h}$ of $h$. However, replacing $\tilde{h}$ with another root
$c\cdot \tilde{h}$ ($c \in \mu_{w'_r}$) changes neither $\alpha^*(\tilde{g}_i)$ nor $\nu_{\bx}(\bof)$.
 \vspace{0.3cm}
 
Let us come back to formula (2.9). Separating the left hand side into the contributions of $\bz \in M(C_{\infty})$ and 
$\bz \in \partial M(C_{\infty})$ and inserting (2.17), we find the ultimate form of the "$k/12$-formula".
 \vspace{0.3cm} 

{\bf 2.19 Theorem.} {\it
Let $\bof = \{f_1,\ldots,f_{r-1}\}$ be an independent family of modular forms for $\Ga$ with weights $k_1,\ldots,k_{r-1}$. For
$\bo \in \OM$ let $s(\bo)$ be the $\gcd$ of $\{i~|~1 \leq i \leq r,\:g_i(\bo) \not=0\}$, and put $w_i = q^{i}-1$ ($i \in \N$).
Then
 $$\underset{\bo \in \Ga\setminus\OM}{\sum} \frac{\nu_{\bo}(\bof)}{w_{s(\bo)}/w_1} + \underset{\bz \in \partial M(C_{\infty})}{\sum}
  \nu_{\bz}(\bof) = \underset{1 \leq i \leq r-1}{\prod} k_i/\underset{2 \leq i \leq r}{\prod} w_i,$$
  where the left hand sum is over a system of representatives for $\Ga\setminus \OM. \:\Box$} 
   \vspace{0.3cm}
   
 {\bf 2.20 Remark} (Comparison with 0.2). Suppose that $r=2$. Then $\bof = \{f\}$ with a single modular form,
 $\nu_{\om}(f)$ is the usual vanishing order of $f$ at $\om \in \OM^2 = C_{\infty}\setminus K_{\infty}$, and 
 $\overline{M} = {\rm Proj}(C_{\infty}[g,\De])$ with $g = g_1$, $\De = g_2$, $\BP = {\rm Proj}(C_{\infty}[X_1,X_2])$, 
 $X_1^{q-1} = g$, $X_2^{q^2-1} = \De$. Hence, for the points ''at infinity'' $\bx = (1:0)\in \BP$, $\bz := \pi(\bx) = (1:0) \in 
 \partial M$, $\nu_{\bx}(\De) = q^2-1$, $\nu_{\bz}(\De)=1$ and correspondingly $\nu_{\bx}(h) = q+1$,
 $\nu_{\bz}(h) = 1/(q-1)$. Each modular form $f$ has a power series expansion 
  $$f(\om) = \underset{i\geq 0}{\sum} a_itî(\om)$$
 with respect to the uniformizer $t(\om)$ ''at infinity" (\cite{4}, Sect.\,5), where $a_i$ vanishes if $f$ is of type zero and 
 $i$ is not divisible by $q-1$. The $\nu_{\infty}(f)$ of (0.2) is the vanishing order with respect to $t$, while the present
 $\nu_{\bz}(f) = (q-1)^{-1} \nu_{\infty}(f)$ is the (virtual, possibly fractional) vanishing order with respect to $s := t^{q-1}$.
 This explains the coefficient $(q-1)^{-1}$ in (0.2). Similar reasoning applies to higher ranks $r \geq 2$ and the vanishing
 order of modular forms along the divisor $\partial M^r$. In this case, there is a uniformizer $t$ along $\partial M^r$, with
 respect to which modular forms may be developed, see \cite{1}.
  \vspace{0.5cm}
  
 \section{Examples.}
We demonstrate how the formula of Theorem 2.19 may be applied. In this section, $\bF$ denotes the ''fundamental domain'' for
$\Ga$ on $\OM$ described in \cite {8}. The rank $r$ is always larger or equal to 2. We start with a very simple case, where
all the ingredients of (2.19) are already known.
 \vspace{0.3cm} 
 
{\bf 3.1 Example.} Let $\bof$ be the family $\{g_1,\ldots ,\hat{g}_j,\ldots,g_r\}$, where $\hat{g}_j$ means $g_j$ deleted. Then the right
hand side of (2.19) evaluates to $(q-1)/(q^j-1)$. 
 \vspace{0.3cm}
 
\underline{Case $j=r$.} For each zero $\bo \in \OM$ of $\bof = \{g_1,\ldots,g_{r-1}\}$, $s(\bo) = r$, i.e., the coefficient of
$\nu_{\bo}(\bof)$ in (2.19) is $(q-1)/(q^r-1)$, too, and there is exactly one such zero $\bo \in \Ga\setminus \OM$, with $\nu_{\bo}(\bof) = 1$, in accordance with the non-vanishing of $\underset{1\leq i,\,j\leq r-1}{\det} (\frac{\partial g_i}{\partial \om_j})$ (\cite{8}, 3.15). Such
a zero is represented in $\bF$ by each element of $\OM^r(\F^{(r)})$, which forms a single orbit under
${\rm GL}(r,\F) \hookrightarrow \Ga$ (\cite{6} 2.7). 
 \vspace{0.3cm}

\underline{Case $j <r$.} There is precisely one zero $\bz \in \overline{M}(C_{\infty})$, in fact $\bz \in \partial M(C_{\infty})$,
with coordinates$(0:\cdots:1:\cdots:0)$ ($z_j=1$), and precisely one $\bx \in \BP(C_{\infty})$ above $\bz$. We have 
$\nu_{\bx}(\bof) = \underset{1\leq i\leq r \atop i \not= j}{\prod} w_i$, and so $\nu_{\bz}(\bof) = (q-1)/(q^j-1)$.  Identifying
$\partial M = \partial M^r$ with $\overline{M}^{r-1}$, we could further restrict the location of $\bz$ as in the case
$j=r$.
 \vspace{0.3cm} 
 
(3.2) In the next (much more interesting) examples, we deal with the forms $E_{q^{î}-1}$, the special Eisenstein series, and the
para-Eisenstein series $\alpha_i$, both of weight $q^{i}-1$ and type zero, see \cite{7} 1.12 or \cite{8} Sect.\,3. These satisfy
various identities and recursions, from which we extract the following consequence.
 \vspace{0.3cm}
 
{\bf 3.3 Proposition.} {\it 
Let $S = \{m,m+1,\ldots,m+r-1\}$ be a set of $r$ consecutive natural numbers and $\bo \in \OM$.
Then at least one of $E_{q^{i}-1}(\bo)$ ($i \in S$) doesn't vanish. The same statement holds true for the $\alpha_i$.}

 \begin{proof} 
Suppose that $E_{q^{i}-1}(\bo) = 0$ for all $i\in S$. The recursion \cite{4} 2.7 with $a = T$, together with the fact
$E_{q^{i}-1} = -\beta_i$ implies that $E_{q^{i}-1} (\bo)$ vanishes for all natural $i \geq m$, which is absurd as the logarithm
function $\log_{\La_{\bo}}(z) = -\underset{i\geq 0}{\sum} E_{q^{i}-1}(\bo)z^{q^{i}}$ has finite convergence radius. The argument
for the $\alpha_i$ is similar: If $\alpha_i(\bo) = 0$ for $i\in S$ then $\alpha_i(\bo) = 0$ for all $i \geq m$ by \cite{8} 3.4\,(ii).
This however conflicts with the fact that $e_{\La_{\bo}}(z) = \sum \alpha_i(\bo)z^{q^{i}}$ is not a polynomial.
 \end{proof}

{\bf 3.4 Corollary.} {\it
Let $j \in S$ and $S' := S \setminus \{j\}$. Then the families $E_{q^{i}-1}$ ($i\in S'$) and $\alpha_i$ ($i\in S'$) are 
independent.}
 
 \begin{proof}
If the set of common zeroes in $\overline{M}$ of the $E_{q^{i}-1}$ ($i\in S'$) is infinite, then it has positive dimension and intersects 
non-trivially with $V(E_{q^j-1})$. This contradicts 3.3. The argument for the $\alpha_i$ is identical.
 \end{proof}
 
 {\bf 3.5 Example.} Let $n$ be a natural number larger or equal to $r$ and $S:= \{i\in \N~|~ n-r < i < n\}$. Then
 $\bof := \{E_{q^{i}-1}~|~i\in S\}$ is independent by (3.4). Further,  $\bof$ has no zero at the boundary $\partial M^r$. 
 For suppose $\bz \in \partial M^r(C_{\infty})$ is such a zero. It corresponds to a Drinfeld module $\phi$ of rank $r'<r$, 
 represented by some $\bo \in \OM^{r'}$. If $r \geq 3$ then (3.3) shows that the $r-1 \geq r'$ consecutive 
 $E_{q^{i}-1}(\bo)$ cannot simultaneously vanish. If $r=2$ then $\phi$ has rank one and is isomorphic with the Carlitz
 module. It is well-known that the $E_{q^{i}-1}$ for the Carlitz module don't vanish, see e.g. \cite{4} 4.3. Now (2.19) reads
  $$\underset{\bo \in \Ga\setminus \OM}{\sum} \frac{\nu_{\bo}(f)}{w_{s(\bo)}/w_1} = \underset{1 \leq i \leq r-1}{\prod}
   (\frac{w_{n-r+i}}{w_{i+1}}). \leqno{(3.5.1)}$$
We are going to locate the zeroes of $\bof$ on the fundamental domain 
 $$\bF = \{\bo \in \OM~|~ \om_r,\om_{r-1},\ldots,\om_1 = 1 \mbox{ is a successive minimum basis of }\La_{\bo} \},$$  
see \cite{8} 1.5 and 1.6. In particular $|\om_r| \geq |\om_{r-1}| \geq \cdots \geq |\om_1| = 1$ for $\bo \in \bF$. In what 
follows, "$\equiv$" means congruence in the ring $O_{C_{\infty}}$ of integers in $C_{\infty}$ modulo its maximal ideal,
and $x \longmapsto \overline{x}$ is the reduction map from $O_{C_{\infty}}$ to its residue class field $\overline{\F}$. 
 \vspace{0.3cm}
 
Suppose that $\bo \in \F_{\boldsymbol 0} =\{\bo \in \bF~|~|\om_1| = \cdots = |\bo_r| = 1\}$. Then 
 \begin{eqnarray*}
 E_{q^{i}-1}(\bo) &=& \underset{\ba \in A^r}{\sum{}{'}} (a_1\bo_1+\cdots + a_r\bo_r)^{1-q^{i}} \equiv  \underset{\ba \in \F^r}{\sum{}{'}}
 (a_1\overline{\om}_1+\cdots + a_r\overline{\bo}_r)^{1-q^{i}}\\
 &=:& E_{q^{i}-1}(\overline{\La}_{\overline{\bo}}),
  \end{eqnarray*}
 where $\overline{\La}_{\overline{\bo}}$ is the $\F$-lattice in $\overline{\F}$ generated by the reductions 
 $\overline{\om}_1,\ldots,\overline{\om}_r$ of the $\om_j$. Now by \cite{5} 1.13, $E_{q^{i}-1}(\overline{\La}_{\overline{\bo}})$
 vanishes for $\overline{\bo} \in \OM^r(\F^{(n)})$ (which is in fact the precise vanishing locus of 
 $E_{q^{i}-1}(\overline{\La}_{\overline{\bo}})$, $i \in S$, see \cite{6} 2.9). The reasoning in \cite{7}, Proposition 6.3
 and Lemma 6.4 shows that the functional determinant 
  $$\det_{n-r< i<n\atop 1 \leq j \leq r-1} (\frac{\partial}{\partial \overline{\om}_j} E_{q^{i}-1} (\overline{\La}_{\overline{\bo}})) 
   \in \overline{\F}$$
 doesn't vanish. Therefore, the multidimensional Hensel lemma implies:
  \vspace{0.3cm}
  
 (3.5.2) Given any $\bo \in \OM^r(\F^{(n)})$, there exists a unique common zero $\bzeta_{\bo} \in \bF_{\boldsymbol 0}$
 of the $E_{q^{i}-1}$ ($i\in S$) with reduction $\overline{\bzeta}_{\bo} = \bo$. For two such, $\bzeta_{\bo}$ and
 $\bzeta_{\bo'}$, we have:
  \begin{eqnarray*}
  \bzeta_{\bo} \mbox{ is $\Ga$-equivalent with } \bzeta_{\bo'} & \Leftrightarrow & \bzeta_{\bo} \mbox{ is $G$-equivalent with }
    \bzeta_{\bo'}\\
     & \Leftrightarrow & \bo \mbox{ is $G$-equivalent with $\bo'$} 
  \end{eqnarray*}
 where $G = {\rm GL}(r,\F) \hookrightarrow \Ga$.
  \vspace{0.3cm}
  
 We will show that the contribution of these zeroes, i.e., of $G$-orbits  on $\OM^r(\F^{(n)})$, to the left hand side of
 (3.5.1) equals the right hand side. So these are all the zeroes of $\bof$ modulo $\Ga$, and $\nu_{\bzeta}(\bof) = 1$ for all
 $\bzeta = \bzeta_{\bo}$.
  \vspace{0.3cm}
  
For each divisor $s$ of $g:= \gcd(r,n)$ let $\OM^r(\F^{(n)})(s)$ be the $G$-stable set of $\bo \in \OM^r(\F^{(n)})$ such that
$\# G_{\bo} = w_s = q^s-1$. Its elements correspond to $\F$-lattices which are in fact $\F^{(s)}$-vector spaces. We have
 $$\underset{s|g}{\sum} \#\OM^r(\F^{(n)})(s) = \#\OM^r(\F^{(n)}).\leqno{(3.5.4)}$$
The number of $G$-orbits on $\OM^r(\F^{(n)})(s)$ is $\frac{\#\OM^r(\F^{(n)})(s) w_s}{\#(G)}$. The contribution of all
$G$-orbits on $\OM^r(\F^{(n)})$ to the left hand side of (3.5.1) is therefore
 $$ \underset{s|g}{\sum} \frac{\OM^r(\F^{(n)})(s)\cdot w_s}{\#(G)} (\frac{w_1}{w_s}) = \frac{q-1}{\#(G)} \#\OM^r(\F^{(n)}),
  \leqno{(3.5.5)}$$
 which by an easy calculation agrees with the right hand side of (3.5.1). We have therefore shown the following result.
  \vspace{0.3cm}
  
 {\bf 3.6 Proposition.} {\it  
Let $\bof$ be the family of special Eisenstein series $E_{q^{i}-1}$, where $n-r< i< n$ for some $n\geq r$. Then 
$\bof$ has no common zeroes on $\partial M^r$. The common zeroes of $\bof$ in the fundamental domain $\bF$ lie in
$\bF_{\boldsymbol{0}}= \{\bo \in \bF~|~|\om_1| = \cdots = |\bo_r| = 1 \}$ and correspond via reduction bijectively 
to the set $\OM^r(\F^{(n)})$. Two such zeroes are $\Ga$-equivalent if and only if the corresponding elements of
$\OM^r(\F^{(n)})$ are ${\rm GL}(r,\F)$-equivalent. All these zeroes $\bzeta$ have multiplicity $1$, i.e., 
\mbox{$\nu_{\bzeta}(\bof) = 1$}.}
 \vspace{0.3cm}
 
We refrain from writing down the precise number of such zeroes $\bzeta$ (or of zeroes $\bzeta$ with $s(\bzeta)$ a given divisor
$s$ of $\gcd(r,n)$), which is an easy exercise in manipulating the M\"obius function. 
 \vspace{0.3cm}
 
{\bf 3.7 Example.} Consider the family $\bof = \{\alpha_2,\alpha_3, \cdots,\alpha_r\}$. It is independent without common zeroes
at $\partial M^r$. The right hand side of (2.19) yields 1. By \cite{8}, Theorem 4.8, a zero $\bo$ of $\bof$ in $\bF$ must
satisfy $|\om_1| = |\om_2| = \cdots = |\om_{r-1}| = q$, $\om_1 = 1$, which excludes $s(\bo) > 1$. Therefore, (2.19) implies that
up to $\Ga$-equivalence, there is exactly one zero $\bo$ of $\bof$, necessarily with $\nu_{\bo}(\bof) = 1$, and represented
by an element of 
 $$\bF_{\boldsymbol{k}}= \{\bo \in \bF~|~|\om_1| = |\om_2| = \cdots = |\om_{r-1}|=q,\,\om_1=1\},$$
 $\boldsymbol{k} = (1,1,\ldots,1,0) \in \Z^r$. (The notation is in accordance with \cite{8}.) Such a zero is uniquely determined
 up to the action of the stabilizer

  $$
  \Ga_{\boldsymbol{k}} =\{ 
   \begin{array}{|c|c|} \hline
     \alpha  &\beta \\ \hline 
     0 & \delta  \\ \hline
      \end{array} ~|~\alpha \in {\rm GL}(r-1,\F),\: \beta \in A^{r-1} \mbox{ with entries of degree } \leq 1,\,
       \delta \in \F^*\}
      $$
of $\bF_{\boldsymbol{k}}$ in $\Ga$.
 \vspace{0.3cm}
 
{\bf 3.8 Concluding remark.} The vanishing sets of the Eisenstein series $E_{q^{i}-1}$ and the para-Eisenstein series $\alpha_i$ and
their intersections are $K$-subvarieties of $\overline{M}^r$. Presumably they may be used for the construction of towers of
curves over finite fields with many rational points in the style of \cite{6}, but they are also interesting in their own right.
A first concrete problem will be to find an analogue of Proposition 3.6 for the para-Eisenstein series $\alpha_i$.

  \vspace{0.5cm}
  
  Ernst-Ulrich Gekeler\\
  Fachrichtung Mathematik\\
  Univesit\"{a}t des Saarlandes\\
  Campus E2 4\\
  66123 Saarbr\"{u}cken\\
  Germany\\
  gekeler@math.uni-sb.de


\begin{thebibliography}{99}
  \bibitem[1]{1} Basson, Dirk and Breuer, Florian: On certain Drinfeld modular forms of higher rank. To appear in J.
 Th\'eorie des Nombres de Bordeaux. 
 \bibitem[2]{2} Bourbaki, Nicolas: Alg\`{e}bre commutative, Chapitres 8 et 9, Masson, Paris 1983. 
 \bibitem[3]{3} Eisenbud, David: Commutative Algebra with a view Toward Algebraic Geometry. Graduate Texts in
 Mathematics 150, Springer-Verlag 1995.
  \bibitem[4]{4} Ernst-Ulrich Gekeler: On the coefficients of Drinfeld modular forms. Invent. Math. 93 (1988) 667--700.
  \bibitem[5]{5} Ernst-Ulrich Gekeler: Finite modular forms. Finite Fields Appl. 7 (2001), 553--572.
  \bibitem[6]{6} Ernst-Ulrich Gekeler: Towers of ${\rm GL}(r)$-type of modular curves. To appear in J. reine angew. Math.
  DOI 10.1515/crelle-2007-0012.
  \bibitem[7]{7} Ernst-Ulrich Gekeler: On Drinfeld modular forms of higher rank. To appear in J. Th\'eorie des Nombres de
  Bordeaux.
  \bibitem[8]{8} Ernst-Ulrich Gekeler: On Drinfeld modular forms of higher rank II, arXiv: 1708.04197 (2017).
  \bibitem[9]{9} Hartshorne, Robin: Algebraic Geometry. Graduate Texts in Mathematics 52, Springer-Verlag 1977.
  \bibitem[10]{10} Serre, Jean-Pierre: Cours d'Arithm\'etique. Presses Universitaires de France, Paris 1970.
  \end{thebibliography}
\end{document}